\documentclass[11pt]{amsart}

\def\pr{{\mathbb P}}

\def\Aff{{\mathbb A}}

\def\skp{\hspace{1pt}}

\def\ph{\varphi}

\def\Id{\mathop{\rm Id\skp}}

\def\rk{{\mathop{\rm rk}\skp}}

\def\prj{{\mathop{\rm pr}\skp}}
\def\img{{\mathop{\rm im}\skp}}

\def\ind{{\mathop{\rm ind}\skp}}

\def\Hom{{\mathop{\rm Hom}\skp}}

\def\Spec{{\mathop{\rm Spec}\skp}}
\def\Proj{{\mathop{\rm Proj}\skp}}

\newcommand {\shExt}  {\mathcal{E} \!\textrm{\textit{xt}}}
\newcommand {\shHom}  {\mathcal{H} \!\textrm{\textit{om}}}

\newcommand{\llongrightarrow}{{\begin{picture}(25,6)(2.5,-4)
  \unitlength 1pt\put(0,0){\vector(1,0){30}}\end{picture}}}
\newcommand{\longdownarrow}{{\begin{picture}(0,0)
  \unitlength 1pt\put(0,16){\vector(0,-1){30}}\end{picture}}}
\newcommand{\longuparrow}{{\begin{picture}(0,0)
  \unitlength 1pt\put(0,-16){\vector(0,1){30}}\end{picture}}}
\newcommand{\arrowrtlb}{{\begin{picture}(0,0)
  \unitlength 1pt\put(15,15){\vector(-1,-1){30}}\end{picture}}}
\newcommand{\arrowlbrt}{{\begin{picture}(0,0)
  \unitlength 1pt\put(-15,-15){\vector(1,1){30}}\end{picture}}}
\newcommand{\arrowltrb}{{\begin{picture}(0,0)
  \unitlength 1pt\put(-15,15){\vector(1,-1){30}}\end{picture}}}
\newcommand{\arrowrblt}{{\begin{picture}(0,0)
  \unitlength 1pt\put(15,-15){\vector(-1,1){30}}\end{picture}}}

\newcommand{\diagl}[1]%
  {\makebox[0cm]{${\scriptstyle#1\ }\longdownarrow
  \phantom{\scriptstyle#1\ }$}}
\newcommand{\diagr}[1]%
  {\makebox[0cm]{$\phantom{\ \scriptstyle#1}
  \longdownarrow{\ \scriptstyle#1}$}}
\newcommand{\updiagl}[1]%
  {\makebox[0cm]{${\scriptstyle#1\ }\longuparrow
  \phantom{\scriptstyle#1\ }$}}
\newcommand{\updiagr}[1]%
  {\makebox[0cm]{$\phantom{\ \scriptstyle#1}\longuparrow
  {\ \scriptstyle#1}$}}
\newcommand{\rtlbdiagl}[1]%
  {\makebox[0pt]{${\scriptstyle#1\atop\ \ }\arrowrtlb
  \phantom{\scriptstyle#1\atop\ \ }$}}
\newcommand{\lbrtdiagl}[1]%
  {\makebox[0pt]{${\scriptstyle#1\atop\ \ }\arrowlbrt
  \phantom{\scriptstyle#1\atop\ \ }$}}
\newcommand{\ltrbdiagl}[1]%
  {\makebox[0pt]{${\ \ \atop\scriptstyle#1\ }\arrowltrb
  \phantom{\ \ \atop\scriptstyle#1\ }$}}
\newcommand{\rbltdiagl}[1]%
  {\makebox[0pt]{${\ \ \atop\scriptstyle#1\ }\arrowrblt
  \phantom{\ \ \atop\scriptstyle#1\ }$}}
\newcommand{\rtlbdiagr}[1]%
  {\makebox[0pt]{$\phantom{\ \ \atop\ \scriptstyle#1}
  \arrowrtlb{\ \ \atop\ \scriptstyle#1}$}}
\newcommand{\lbrtdiagr}[1]%
  {\makebox[0pt]{$\phantom{\ \ \atop\ \scriptstyle#1}
  \arrowlbrt{\ \ \atop\ \scriptstyle#1}$}}
\newcommand{\ltrbdiagr}[1]%
  {\makebox[0pt]{$\phantom{\scriptstyle#1\atop\ \ }
  \arrowltrb{\scriptstyle#1\atop\ \ }$}}
\newcommand{\rbltdiagr}[1]%
  {\makebox[0pt]{$\phantom{\scriptstyle#1\atop\ \ }
  \arrowrblt{\scriptstyle#1\atop\ \ }$}}

\newcommand{\calc}{{\mathcal C}}

\newcommand{\cale}{{\mathcal E}}
\newcommand{\calf}{{\mathcal F}}

\newcommand{\cali}{{\mathcal I}}

\newcommand{\calk}{{\mathcal K}}
\newcommand{\call}{{\mathcal L}}

\newcommand{\calo}{{\mathcal O}}

\newcommand{\calq}{{\mathcal Q}}

\newcommand{\calu}{{\mathcal U}}

\newcommand{\calx}{{\mathcal X}}

\newtheorem{prop}{Proposition}[section]
\newtheorem{theorem}[prop]{Theorem}
\newtheorem{lemma}[prop]{Lemma}

\newtheorem{rem}[prop]{Remark}
\newtheorem{expl}[prop]{Example}
\newtheorem{defi}[prop]{Definition}
\newcommand{\pf}{{\em Proof. }}

\newcommand{\Lin}{{\mathop{\rm Lin}}}
\newcommand{\Coh}{{\mathop{\rm Coh}}}

\title[Virtual fundamental classes]{Virtual fundamental classes, global normal cones
         and  Fulton's canonical classes}
\author{Bernd Siebert}
\address{Mathematisches Institut, Universit\"at Freiburg, D-79102 Freiburg}
\email{bernd.siebert@math.uni-freiburg.de}
\date{February 5, 1997; revision of August 14, 2005}
\begin{document}
\maketitle

\tableofcontents
\vspace{4ex}


\noindent
{\Large\bf Introduction}
\vspace{1.5ex}

\noindent
This note, written in January 1997, grew out of an attempt to
understand references \cite{behrend}, \cite{behrendfantechi} and
\cite{litian}. In these papers two related but different methods are
presented for the construction of a certain Chow class on moduli
spaces of stable (parametrized) curves in a projective manifold $V$,
called virtual fundamental class. This class replaces the usual
fundamental class of these spaces in the definition of basic
enumerative invariants of $V$ involving curves, called Gromov-Witten
(GW-) invariants. They are invariant under smooth deformations of
$V$.

Both approaches are based on a globalization of the concept of
normal cones of germs of the space under study inside some
modelling space, that is $C_{U|M}$ for $U\subset X$ open with $\iota: U
\hookrightarrow M$ and $M$ smooth over $k$. The essential idea of
using bundles of cones inside a vector bundle for globalizing virtual
fundamental classes is due to Li and Tian. The data needed to glue
differs however somewhat in the two constructions.

A proper understanding of the relationship between the two
approaches seemed necessary for finding the natural framework for
comparison of algebraic virtual fundamental classes with the author's
definition in \cite{si1} of virtual fundamental classes in the
symplectic context \cite{si3}.

In a first step Behrend and Fantechi use a generalization of the
concept of scheme, called Artin stacks, to make sense of the quotient
$C_{U|M}/T_M|_U$. These quotients being unique up to canonical
isomorphism they glue to an Artin (cone) stack ${\calc}_X$
intrinsically associated to any $X$. In a second step they need a
morphism $\ph^\bullet: [\calf^{-1} \rightarrow \calf^0] \rightarrow
\call_X^\bullet$ (in the derived category) from a two-term complex of
locally free sheaves to the cotangent complex, inducing an
isomorphism in $H^0$ and an epimorphism in $H^{-1}$, to cook up an
ordinary cone $C(\ph^\bullet)\subset F_1$, $F_1$ the vector bundle
associated to $\calf^{-1}$. Intersection with the zero section
finally produces the virtual fundamental class. So a priori the
latter depends on the choice of $\ph^\bullet$.

This is a very natural and mature approach, that clearly separates the
globalization process of the normal cone from the construction of the
virtual fundamental class. A possible disadvantage is that in dealing
with Artin stacks some of the necessary verifications become rather
technical, non-geometric in nature.

Li and Tian circumvent the morphism to the cotangent complex by
introducing the notion of ``perfect tangent obstruction complex''.
This is a morphism $\calf^{-1}\rightarrow\calf^0$ of locally free
sheaves on $X$ with kernel and cokernel being tangent and obstruction
spaces for morphisms to $X$, compatible with base change. Using
relative, formal, ``Kuranishi families'' as an intermediate object
they construct a well-defined cone $C\subset F_1$. Another, less
important difference to \cite{behrendfantechi} is the use of an
absolute obstruction theory instead of one relative to the space of
pre-stable curves.

In a previous version of \cite{litian} the slightly stronger claim
was made that already a presentation $\calf^{-1}\rightarrow\calf^0$
of $\Omega_X$ should suffice to construct the cone. In trying to
understand this statement I was lead to the problem of reformulating
\cite{behrendfantechi} from the point of view of gluing local cones.
Since in the latter reference Artin stacks are used only as
book-keeping device rather than as actual spaces it should not come
as surprise that one can get along without them (this has already
been indicated in op.cit.). Contrary to what I expected, things can
be formulated in a rather elegant but direct way via some \emph{yoga
of cones bundles}. This part of the paper (Sections~2 and 3) is just
a down-to-earth reformulation of (parts of) Sections~2--4 of
\cite{behrendfantechi}. Section~1 presents the necessary notations
concerning cones and linear spaces, the latter being a convenient way
of looking at coherent sheaves for our purposes.

In Section~4 we establish a closed formula for virtual fundamental
classes involving only the scheme-theoretic structure of $X$ via
Fulton's canonical class (Definition~\ref{fultons_class}) and the Chern
class of the virtual bundle $F_0-F_1$, $F_i$ the vector bundle
associated to $\calf^{-i}$ (Theorem~\ref{closed_formula}). This
formula was actually found by the author in summer 1995 while
searching for a purely algebraic definition of GW-invariants. It should
be useful for computations, see the author's recent little survey
\cite{si2}.

A few remarks on GW-theory are in order. First, today I consider the
yoga of cone bundles in Sections~2 and 3 as one ingredient for the
most economic path to algebraic Gromov-Witten invariants. The
other ingredients are going over to Deligne-Mumford stacks, and
replacing the morphism to the cotangent complex by an obstruction
theory. If the latter is defined similar to \cite{artin}, 2.6, rather
than in \cite{litian}, one can show \cite{si4} that it is locally
nothing but a morphism to the cotangent complex as in
\cite{behrendfantechi}. Hence the yoga of cone bundles applies to
produce the virtual fundamental class.

Second, I would like to illustrate the perspective of the content of
Section~4 by the following formula for virtual fundamental classes in
Gromov-Witten theory.
\begin{theorem}
  Let $V$ be a projective variety, smooth over a field $K$ of
  characteristic $0$, and $R\in A_1(V)$, the first Chow group. If
  $g=0$ or $\calc:= \calc_{R,g,k}(V)$, the moduli space of stable
  curves $(C,{\bf x},\ph:C\rightarrow V)$ in $V$ of genus $g$ with
  $k$ marked points ${\bf x}= (x_1,\ldots,x_k)$  and $\ph_*[C]=R$, is
  embeddable into a space smooth over ${\mathfrak{M}}_{g,k}$, then
  the virtual fundamental class relevant for GW-invariants is
  \[
    [\![\calc]\!]\ =\ \Big\{c(\ind_{R,g,k}^V)^{-1}\cap
    c_F(\calc/{\mathfrak{M}}_{g,k})\Big\}_{d(V,R,g,k)}\, .
  \]
  Here $\{\,.\,\}_d$ denotes the $d$-dimensional part of a cycle,
  ${\mathfrak{M}}_{g,k}$ is the Artin stack of $k$-pointed pre-stable curves
  of genus $g$, $d(V,R,g,k)= c_1(V)\cdot R +(1-g) \dim V+ 3g-3$ is the
  expected dimension, and $c_F(\calc/ {\mathfrak{M}}_{g,k})$ is Fulton's
  canonical class for $\calc$ relative ${\mathfrak{M}}_{g,k}$.
  \qed
\end{theorem}
Here $\ind_{R,g,k}^V=F_0-F_1$ is the virtual vector bundle associated
to the partial resolution $\varphi^\bullet: [\calf^{-1}\to
\calf^0]\to \call^\bullet_{\calc/ {\mathfrak{M}}_{g,k}}$ mentioned in
the introduction. It represents the (domain of the) perfect relative
obstruction theory $(R\pi_*(f^T_V))^\vee$ of Behrend \cite{behrend}
in $K^0(\calc)$.

A note on categories: To keep things simple we work here in the
category of schemes of finite type over a field $k$, not necessarily
algebraically closed or of characteristic 0. The extension to other
base schemes is straightforward. For the purpose of GW-theory one
also has to replace schemes by (generalizations of) orbifolds, that
is Deligne-Mumford stacks in the algebraic category or analytic
orbispaces in an analytic context. Again, our results can be easily
adapted to these categories. 

For GW-theory this is still not sufficient, because
${\mathfrak{M}}_{g,k}$ is an Artin stack rather than Deligne-Mumford.
One can nevertheless give a construction of the relevant cone without
ever really using Artin stacks. For instance, Fulton's canonical
class relative ${\mathfrak{M}}_{g,k}$ has the following simple
definition: Embed $\calc$ into a smooth Deligne-Mumford $k$-stack
$N$. In the important case $g=0$ one could take $\calc_{\iota_*R,g,k}
(\pr^N)$, if $\iota: V\hookrightarrow \pr^N$ is a closed embedding.
Let $q:\calu \rightarrow \calc$ be the universal curve. Then the
pull-back of the (virtual) tangent bundle of ${\mathfrak{M}}_{g,k}$
is
\[
  T_{\mathfrak{M}} := \shExt^1_q (\omega_{\calu /\calc},\calo_\calu)-
  \shExt^0_q (\omega_{\calu /\calc},\calo_\calu)\,,
\]
as an element of $K^0(\calc)$
(the $\shExt^i_q$ are the derived functors of $q_*\circ\shHom$) and
\begin{eqnarray*}
  c_F(\calc/ {\mathfrak{M}}_{g,k})&=&
  \Big(c(T_{\mathfrak{M}})^{-1} \cup c(T_N)\Big) \cap s(C_{\calc/N})\\
  &=& c(T_{\mathfrak{M}})^{-1}\cap c_F(\calc)\,.
\end{eqnarray*}
Here all sheaves and classes have to be understood in the sense of
Deligne-Mumford stacks. With these remarks understood the theorem
is a special case of Theorem~\ref{closed_formula}.

In preparing this paper discussions with H.\ Flenner and S.\
Schr\"oer have been helpful. I am grateful to the referee for a very
attentive reading of the manuscript and several competent
suggestions.


\section{Cones}
\subsection{Linear spaces}
For any algebraic $k$-scheme $X$, $\Aff^1_X=X\times\Aff^1_k$
has the structure of a {\em ring over} $X$: There are morphisms
\[
  \alpha:\Aff^1_X\times_X\Aff^1_X\longrightarrow\Aff^1_X\, ,\quad
  \iota:\Aff^1_X\longrightarrow\Aff^1_X\, ,\quad
  \mu:\Aff^1_X\times_X\Aff^1_X\longrightarrow\Aff^1_X\, ,
\]
and sections $n,e:X\rightarrow\Aff^1_X$ fulfilling the usual
commutative ring axioms with $\alpha$ as addition, $\iota$ as
additive inverse, $\mu$ as multiplication and $n$, $e$ as neutral
elements for $\alpha$, $\mu$.

A {\em linear space over $X$} is an $\Aff^1_X$-module of finite type
over $X$, that is, an affine morphism $\pi:L\rightarrow X$ of
finite type together with morphisms
\[
  a:L\times_X L\longrightarrow L,\quad m:\Aff^1_X\times_X L
  \longrightarrow L
\]
and a zero section $z:X\rightarrow L$, fulfilling the usual module
axioms relative $X$, that is $m\circ(\mu\times\Id_L) =
m\circ(\Id_{\Aff^1_X}\times m)$ as maps from $\Aff^1_X \times_X
\Aff^1_X\times_X L$ to $L$ etc.\footnote{``Linear space (over $X$)'' or
``linear fiber space'' (``Linearer Faserraum'') seem to be the classical
notation for the ``abelian cones'' of \cite{behrendfantechi}} By abuse
of notation we just write $L$ for the tuple $(\pi,a,m)$. In the sequel
we will restrict to linear spaces that are {\em representable}, that
is, which locally are closed subspaces of vector bundles with induced
linear structure. With the obvious notion of homomorphism of linear
spaces over $X$ we get the category $\Lin(X)$ of representable linear
spaces over $X$.

There is an anti-equivalence of categories
\[
  \Lin(X)\longrightarrow\Coh(X)
\]
to the category of coherent $\calo_X$-modules \cite{EGA-II},\S1.7: On
objects this associates to $L\in \Lin(X)$ the sheaf
$\Hom_{\Lin(X)}(L,\Aff^1_X)$; in the other direction,
$\calf\in\Coh(X)$ corresponds to
\[
  L(\calf)\ :=\ \Spec_{\calo_X}S^\bullet\calf\, ,
\]
where $S^\bullet\calf$ is the symmetric algebra over the
$\calo_X$-module $\calf$. For example, the addition operation $a:
L(\calf)\times_X L(\calf)\to L(\calf)$ comes from the diagonal
morphism $\calf\to \calf\oplus \calf$, $f\mapsto (f,f)$ by
application of the functor $L=\Spec_{\calo_X}\circ S^\bullet$. Note
that a vector bundle $E$ corresponds to the locally free sheaf
$\calo(E^\vee)$, $E^\vee$ the {\em dual} bundle.

Representable linear spaces are thus just another way to look at
coherent sheaves. We will jump freely between both descriptions and
use whichever seems more appropriate in a particular context. Note
also that $\Lin(X)$ is an abelian category, so it makes sense to talk
about monomorphisms, epimorphisms and exact sequences. A
monomorphism $\Phi:E\rightarrow F$ of linear spaces corresponds to
an epimorphism $\ph:\calf \rightarrow\cale$ of sheaves and is thus a
closed embedding of schemes. An epimorphism $\psi:F\rightarrow G$ of
linear spaces, however, need not be a surjection of schemes (consider
the inclusion $\psi:\cali\hookrightarrow\calo_X$ for any nontrivial
ideal sheaf $\cali$).


\subsection{Cones}
A {\em cone} $C$ over $X$ is a scheme of the form $\Spec_{\calo_X}
S^\bullet$ where $S^\bullet= \oplus_{d\ge0}S^d$ is a graded
$\calo_X$-module with $S^0=\calo_X$ and $S^\bullet$ generated by
$S^1\in\Coh(X)$. $S^\bullet$ as graded algebra is not in general
determined up to isomorphism by the scheme $C$ over $X$. For the
grading one needs to distinguish the generating submodule
$S^1=\calf$, or, equivalently, a closed embedding $C\hookrightarrow
L(\calf)$ into a linear space. Such datum could be called {\em
polarization} of $C$. We will only deal with polarized cones in the
sequel.
\begin{expl}\label{normal_cone}\em
If $X$ is a closed subscheme of an algebraic $k$-scheme
$M$ with ideal sheaf $\cali$ then the cone
\[
  C_{X|M}=\Spec_{\calo_X}\big(\oplus_{d\ge0}\cali^d/\cali^{d+1}\big)
\]
over $X$ is called {\em normal cone} to $X$ in $M$. $C_{X|M}$ is
naturally embedded into the {\em normal space}
$N_{X|M}=L(\cali/\cali^2)$ of $X$ in $M$ (to avoid confusion with
$(\cali/\cali^2)^\vee$, I would rather not call $N_{X|M}$ normal
sheaf as in \cite{behrendfantechi}).
\qed
\end{expl}
If $C$, $C'$ are cones over $X$, then so is $C\oplus C':=C\times_X C'$.

To a polarized cone $C=\Spec S^\bullet$ is associated a Chow class on
$X$, its Segre class
\[
  s(C)\ :=\ \sum_{r\ge0}p_*\big(\xi^r\cap[\pr(C)]\big),
\]
where $p:\pr(C) :=\Proj S^\bullet\rightarrow X$ is the projection and
$\xi = c_1(\calo_{\pr(C)}(1))$.

We propose the following formulation of the concept of exact
sequence of cones \cite[Expl.4.1.6]{fulton}.
\begin{defi}\rm
  Let
  \[
    \mbox{\hspace*{4cm}}
    0\longrightarrow E\stackrel{\Phi}{\longrightarrow}F
    \stackrel{\Psi}{\longrightarrow}Q\longrightarrow 0
    \mbox{\hspace*{4cm}}(*)
  \]
  be an exact sequence of linear spaces. Let $C\subset Q$ be a cone
  and set $\tilde C:=\Psi^{-1}(C)$. Then $(*)$ restricts to
  \[
    0\longrightarrow E\longrightarrow\tilde C\longrightarrow C
    \longrightarrow 0 \quad.
  \]
  Sequences of cones of this form will be called {\em exact}.
\qed
\end{defi}
\begin{rem}\em
Exact sequences of cones might not be very useful unless $(*)$
splits locally. In this case $\tilde C$ is locally of the form $C\oplus E$,
and as in \cite[Expl.4.1.6]{fulton} one can show $s(\tilde C)=s(C\oplus
E)$. In the non-split case a convenient way to relate the Segre classes of
$\tilde C$ and $C$ seem to be unknown.

But note that if $E$ is a vector bundle $(*)$ always splits locally,
and so we retrieve the definition of exact sequences of cones as in
\cite{fulton}.
\qed
\end{rem}

For an exact sequence of cones as in the definition $\tilde C$ is
preserved by the additive action of $E$ on $F$. In other words,
$\tilde C$ wears the structure of an $E$-module. More generally, if
$\Phi:E\rightarrow F$ is a homomorphism of linear spaces and
$C\subset F$ is a  cone then $C$ is called {\em $E$-cone} if $C$ is
an $E$-module via $\Phi$, that is if $C$ is preserved by the additive
action of $E$ on $F$ induced by $\Phi$.
\begin{expl}\rm\label{C_Z|X}
In the situation of Example~\ref{normal_cone} $C_{X|M}$ is a
$T_M|_X$-cone via the natural homomorphism $\Phi:T_M|_X
\rightarrow N_{X|M}$. On the sheaf level this action of $T_M|_X$ is
\begin{eqnarray*}
\bigoplus_d \cali^d/\cali^{d+1}&\longrightarrow&
S^\bullet \Omega_M|_X \otimes \bigoplus_d \cali^d/\cali^{d+1},
\end{eqnarray*}
where for $f_i\in \cali/\cali^2$ the image of $f_1\cdot\ldots\cdot f_d$ in
the direct summand $S^e \Omega_M|_X \otimes \cali^{d-e}/\cali^{d-e+1}$
of the target is the sum over all partitions $\{i_1,\dots,i_e\}$,
$\{j_1,\dots,j_{d-e}\}$ of $\{1,\dots,n\}$ of terms
\begin{eqnarray*}
\textrm{d}f_{i_1}\cdot \ldots\cdot \textrm{d}f_{i_e}\otimes
f_{j_1}\cdot\ldots\cdot f_{j_{d-e}}.
\end{eqnarray*}\vspace{-6ex}

\qed
\end{expl}
There are many examples of morphisms of linear spaces $E\rightarrow F$ and
$E$-cones $C\subset F$ that do not descend to the quotient $Q=E/F$, 
for instance the examples in Remark~\ref{coh_cones},3 and in
Remark~\ref{C_X|M_descend?}. However, there is one important class
of morphisms where it is always possible, namely for locally split
monomorphisms. We first treat the split case:
\begin{lemma}\label{C_in_EplusF}
  Let $\cale$, $\calf\in\Coh(X)$ and $E=L(\cale)$, $F=L(\calf)$ the 
  corresponding linear spaces over $X$ and $C\subset E\oplus F$ an
  $E$-invariant closed subscheme with respect to the action of
  $E$ on the first summand.

  Then $C$ is of the form $E\oplus\bar C$ for some uniquely
  determined closed subscheme $\bar C\subset F$.
\end{lemma}
\pf
The statement is local in $X$, so we may assume $X=\Spec A$, 
$E=\Spec A[{\bf X}]/\langle{\bf e}\rangle$, $F=\Spec A[{\bf Y}]/
\langle{\bf f}\rangle$ with ${\bf X}=(X_1,\ldots X_r)$, ${\bf Y}=
(Y_1,\ldots Y_s)$ and ${\bf e}=(e_1,\ldots,e_k)$, ${\bf
f}=(f_1,\ldots,f_l)$ tuples of linear forms with coefficients in $A$,
$\langle{\bf e}\rangle$, $\langle{\bf f}\rangle$ the ideals
generated by their entries. Then $C=\Spec A[{\bf X},{\bf Y}]/I$ with $I$
an ideal containing $\langle{\bf e}\rangle+\langle{\bf f}\rangle$.

The only possible candidate for $\bar C$ is the intersection of $C$
with $0\oplus F$, that is $\bar C=\Spec A[{\bf X},{\bf Y}]/(I+\langle{\bf
X}\rangle)=\Spec A[{\bf Y}]/\bar I$, with $\bar I=\{f(0,{\bf Y})\mid f({\bf
X},{\bf Y})\in I\}$. We have to show that $I=\langle\,\bar
I\,\rangle+\langle{\bf e}\rangle$.

$C$ to be $E$-invariant means that for any $\displaystyle f({\bf X},
{\bf Y})=\sum_{M,N}a_{MN}{\bf X}^M{\bf Y}^N\in I$
\[
  \mbox{\hspace*{2cm}}
  f({\bf X}+{\bf X'},{\bf Y})\ =\ \sum_{M,N}a_{MN}({\bf X}+{\bf X'})^M
  {\bf Y}^N \ \in\ \langle\,I\,\rangle+\langle\ph({\bf e})
  \rangle\mbox{\hspace*{2cm}}(*)
\]
holds in $A[{\bf X},{\bf X'},{\bf Y}]$, where $\ph:A[{\bf X}] \rightarrow
A[{\bf X'}]$, $X_\mu\mapsto X_\mu'$. Modulo $\langle{\bf X}\rangle$
this says
\[
  f({\bf X'},{\bf Y})\in\langle\,\bar I\,\rangle+\langle\ph({\bf e})\rangle
\]
in $A[{\bf X'},{\bf Y}]$. Replacing $\bf X'$ by $\bf X$ we thus get
$I\subset\langle\,\bar I\,\rangle+\langle{\bf e}\rangle$.

For the other direction we look at $(*)$ modulo ${\bf X}+{\bf X'}$
to conclude
\[
  f(0,{\bf Y})\ =\ \sum_N a_{0N}{\bf Y}^N
  \ \in\ I+{\bf e}\ =\ I
\]
for any $f\in I$, that is $\bar I\subset I$.
\qed

\begin{prop}\label{quot_cone}
  Let
  \[
    \mbox{\hspace*{5cm}}0\longrightarrow F\longrightarrow E
    \stackrel{q}{\longrightarrow}Q
    \mbox{\hspace*{5cm}}(*)
  \]
  be an exact sequence of linear spaces with $F$ a vector bundle,
  and let $C\subset E$ be an $F$-cone.

  Then there exists a unique cone $\bar C\subset Q$ such that
  $(*)$ induces an exact sequence of cones
  \[
    0\longrightarrow F\longrightarrow C\longrightarrow\bar C
    \longrightarrow0\, .
  \]
  In particular, $C$ descends to $Q$: $C=q^{-1}(\bar C)$.
\end{prop}
\pf
By replacing $Q$ by the closed subspace $E/F\subset Q$ we
may assume $q$ to be an epimorphism. Then, since $F$ is a vector
bundle, locally $(*)$ splits and we may apply the previous lemma to
construct $\bar C\subset Q$.
\qed
\smallskip

\noindent
In other words, the proposition says that $\bar C$ is the
scheme-theoretic quotient of $C$ by the free action of $F$. This is a
convenient way to think about $\bar C$.


\section{Going up and down for $E$-cones}
In this section we investigate the behavior of $E$-cones
under morphisms of two-term complexes, that is commutative squares,
in $\Lin(X)$. If $\Phi_\bullet=(\Phi_0,\Phi_1):
F_\bullet=(F_0\rightarrow F_1)\rightarrow (E_0\rightarrow E_1)$ is
such a morphism the corresponding morphism of coherent sheaves will
be written $\ph^\bullet=(\ph^{-1},\ph^0):\cale^\bullet
=(\cale^{-1}\rightarrow\cale^0)\rightarrow\calf^\bullet=
(\calf^{-1}\rightarrow\calf^0)$. Then $\Phi_i=L(\ph^{-i})$,
$E_i=L(\cale^{-i})$, $F_i=L(\calf^{-i})$ for $i=0,1$.


\subsection{Going up}
\begin{lemma}
  Let $\Phi_\bullet:F_\bullet\rightarrow E_\bullet$ be a commutative 
  square in $\Lin(X)$, and $C\hookrightarrow E_1$ an $E_0$-cone. Then
  $\Phi_1^{-1}(C)\hookrightarrow F_1$ is an $F_0$-cone.
\end{lemma}
\pf
Consider the diagram
\[
\begin{array}{ccc}
F_0\oplus F_1&\stackrel{\alpha}{\llongrightarrow}&F_1\\[10pt]
\diagl{\Phi_0\oplus\Phi_1}&&\diagr{\Phi_1}\\[10pt]
E_0\oplus E_1&\stackrel{\alpha'}{\llongrightarrow}&E_1
\end{array}
\]
with horizontal arrows the morphisms defining the $F_0$- and 
$E_0$-module structures on $F_1$ and $E_1$ respectively. By
hypothesis $E_0\oplus C$ is a closed subscheme of
$(\alpha')^{-1}(C)$. Thus $F_0\oplus\Phi_1^{-1}(C)
=(\Phi_0\oplus\Phi_1)^{-1} (E_0\oplus C)$ is a closed subscheme of
$\alpha^{-1} (\Phi_1^{-1}(C))$.
\qed

By this lemma we are able to make the following definition.
\begin{defi}\rm
(going up) Let $\Phi_\bullet:F_\bullet\rightarrow E_\bullet$ be a 
commutative square in $\Lin(X)$ and $C\subset E_1$ an $E_0$-cone.
Then the $F_0$-cone
\[
  \Phi_\bullet^!(C):=\Phi_1^{-1}(C)
\]
in $F_1$ is called {\em pull-back} of $C$ under $\Phi_\bullet$.
\qed
\end{defi}
The pull-back depends only on the homotopy class of $\ph^\bullet$
(or $\Phi_\bullet$).
\begin{prop}
  Let $\ph^\bullet, \psi^\bullet:\cale^\bullet=[\cale^{-1} \stackrel{d}{\rightarrow}
  \cale^0] \rightarrow \calf^\bullet$ be homotopic commutative squares in
  $\Coh(X)$. Then for any $E_0$-cone $C\subset E_1$
  \[
    \Phi_\bullet^!(C)\ =\ \Psi_\bullet^!(C)\, .
  \]
\end{prop}
\pf
Let $k:\cale^0\to \calf^{-1}$ be a homotopy: $\psi^{-1} =
\ph^{-1}+k\circ d$, $\psi^0=\ph^0+ d\circ k$.
Writing $K=L(k)$ and $\alpha:E_0\oplus E_1\rightarrow E_1$ for the 
structure map, $\Psi_1$ may be decomposed into
\[
  F_1\stackrel{(K,\Phi_1)}{\longrightarrow}E_0\oplus E_1
  \stackrel{\alpha}{\longrightarrow}E_1\, .
\]
Since $E_0\oplus C\subset\alpha^{-1}(C)$, $(K,\Phi_1)^{-1}
(E_0\oplus C) =\Phi_1^{-1}(C)$ is a closed subscheme of
$\Psi_1^{-1}(C)$. But the claim is symmetric in $\Phi_\bullet$,
$\Psi_\bullet$, hence $\Phi_1^{-1}(C)=\Psi_1^{-1}(C)$.
\qed

The next result about functoriality of going up follows directly from
the definition.
\begin{prop}\label{gu_functorial}
  Let $\Phi_\bullet:E_\bullet\rightarrow F_\bullet$, $\Psi_\bullet:
  F_\bullet\rightarrow G_\bullet$ be commutative squares of linear
  spaces and $C\subset G_1$ a $G_0$-cone. Then
  \[
    (\Psi_\bullet\circ\Phi_\bullet)^!(C)\ =\ \Phi_\bullet^!\circ
    \Psi_\bullet^!(C)\, .
  \]
  \vspace{-8ex}
  
  \qed
\end{prop}


\subsection{Going down}
Going down, or push-forward, of $F_0$-cones in $F_1$ to $E_1$ is a 
little more subtle. The central tool will be
Proposition~\ref{quot_cone}. To make this proposition applicable we
need a little lemma.
\begin{lemma}\label{ass_exact_sequence}
  Let $\ph^\bullet:(\cale^{-1}\stackrel{d}{\rightarrow}\cale^0)
  \rightarrow (\calf^{-1}\stackrel{d'}{\rightarrow}\calf^0)$ be a
  commutative square in $\Coh(X)$. Then the complex
  \[
    0\longrightarrow\cale^{-1}\stackrel{(d,\ph^{-1})}{\longrightarrow}
    \cale^0\oplus\calf^{-1}\stackrel{\ph^0\circ\prj_1
    -d'\circ\prj_2}{\longrightarrow}\calf^0\longrightarrow0
  \]
  is exact at
  \vspace{1ex}

  \begin{tabular}[t]{lll}
  i) & $\calf^0$ & iff $H^0(\ph^\bullet)$ is surjective\\
  ii) & $\cale^0\oplus\calf^{-1}$ & iff $H^0(\ph^\bullet)$ is injective
      and $H^{-1}(\ph^\bullet)$ is surjective\\
  iii) & $\cale^{-1}$ & iff $H^{-1}(\ph^\bullet)$ is injective.
  \end{tabular}
\end{lemma}
\pf
Chase the diagram
\[
\begin{array}{rcccccccccccl}
\phantom{.}&
0&\llongrightarrow&\calk&\llongrightarrow&\cale^{-1}
  &\stackrel{d}{\llongrightarrow}&\cale^0
  &\llongrightarrow&\calq&\llongrightarrow&0\\[8pt]
&&&\diagl{H^{-1}(\ph^\bullet)}&&\diagl{\ph^{-1}}&&
  \diagr{\ph^0}&&\diagr{H^0(\ph^\bullet)}\\[12pt]
&0&\llongrightarrow&\calk'&\llongrightarrow&\calf^{-1}
  &\stackrel{d'}{\llongrightarrow}&\calf^0
  &\llongrightarrow&\calq'&\llongrightarrow&0&.
\end{array}
\]
\vspace{-6ex}

\qed
\vspace{2ex}

\noindent
If $\ph^\bullet$ is a quasi-isomorphism we thus get exactness of the
stated complex. And $\ph^\bullet$, viewed as a commutative square, is
cartesian ($\cale^{-1}=\cale^0\oplus_{\calf^0}\calf^{-1}$) iff
$H^0(\ph^\bullet)$ is injective and $H^{-1}(\ph^\bullet)$ is an
isomorphism, and it is cocartesian ($\calf^0=(\cale^0\oplus
\calf^{-1})/ \cale^{-1}$) iff $H^0(\ph^\bullet)$ is an isomorphism
and $H^{-1}(\ph^\bullet)$ is surjective. Assume now that $F_0$ is a
vector bundle and that $\Phi_\bullet:[F_0
\stackrel{D'}{\rightarrow}F_1]\rightarrow[E_0
\stackrel{D}{\rightarrow}E_1]$ induces an isomorphism on $H^0$ and a
closed embedding of linear spaces on $H^1$. If these conditions are
satisfied we say that {\em going down is applicable to
$\Phi_\bullet$}. Then
\[
  0\longrightarrow F_0\stackrel{(\Phi_0,-D')}{\longrightarrow}
  E_0\oplus F_1 \stackrel{q}{\longrightarrow}E_1
\]
is exact (Lemma~\ref{ass_exact_sequence}, $q=D\circ\prj_1
+\Phi_1\circ\prj_2$) and we may apply Proposition~\ref{quot_cone}.
\begin{defi}\label{going_down}\rm
  (going down) Let $\Phi_\bullet:F_\bullet\rightarrow E_\bullet$ be
  a  commutative square in $\Lin(X)$, to which going down is
  applicable (see above), and let $C\subset F_1$ be an $F_0$-cone.
  The unique cone $\bar C\subset\img q\subset E_1$ with $q^{-1}(\bar
  C)=E_0\oplus C$, which exists by Proposition~\ref{quot_cone}, is
  called {\em push-forward} of $C$ by $\Phi_\bullet$, denoted
  $(\Phi_\bullet)_!(C)$.
\qed
\end{defi}

Note that $(\Phi_\bullet)_!(C)$ is actually an $E_0$-cone because 
$E_0\oplus C$ is one. And by Proposition~\ref{quot_cone}:
\begin{prop}\label{ex_seq_cones}
If going down is applicable to $\Phi_\bullet:F_\bullet\rightarrow 
E_\bullet$, and $C\subset F_1$ is an $F_0$-cone, there is an exact
sequence of cones
\[
  0\longrightarrow F_0\longrightarrow E_0\oplus C\longrightarrow
  (\Phi_\bullet)_!(C)\longrightarrow0\, .
\]
\vspace{-8ex}

\qed
\end{prop}
\begin{rem}\rm
  Local freeness of $F_0$ (or local splittability of the relevant exact 
  sequence of linear spaces) seems to be indispensable, since otherwise
  $E_0\oplus C$ need not descend to $E_1$. See
  Remark~\ref{coh_cones},3 for a related example.
\qed
\end{rem}
As with going up, going down depends only on the homotopy class of 
$\Phi_\bullet$.
\begin{prop}
  Let $\Phi_\bullet$, $\Psi_\bullet:F_\bullet\rightarrow E_\bullet$ be 
  homotopic morphisms of commutative squares in $\Lin(X)$ and
  $C\subset F_1$ an $F_0$-cone. If going down is applicable to
  $\Phi_\bullet$ (or, equivalently, to $\Psi_\bullet$) then
  \[
    (\Phi_\bullet)_!(C)\ =\ (\Psi_\bullet)_!(C)\, .
  \]
\end{prop}
\pf
Let $K:F_1\rightarrow E_0$ be a homotopy between $\Phi_\bullet$ 
and $\Psi_\bullet$, that is $\Psi_0=\Phi_0+K\circ D'$, $\Psi_1=
\Phi_1+D\circ K$ ($D:E_0\rightarrow E_1$, $D':F_0 \rightarrow F_1$ the
differentials). Then the following diagram
\[
\begin{array}{ccccccc}
0&\llongrightarrow&F_0&\stackrel{\makebox[0pt]{$\scriptstyle
  (\Phi_0,-D')$}}{\llongrightarrow}&E_0\oplus F_1
  &\stackrel{q_\Phi}{\llongrightarrow}&E_1\\[10pt]
  &&\diagl{\Id}&&\diagl{\chi}&&\diagl{\Id}\\[10pt]
0&\llongrightarrow&F_0&\stackrel{\makebox[0pt]{$\scriptstyle
  (\Psi_0,-D')$}}{\llongrightarrow}&E_0\oplus F_1
  &\stackrel{q_\Psi}{\llongrightarrow}&E_1
\end{array}
\]
with $\chi=(\prj_1-K\circ\prj_2,\prj_2)$, $q_\Phi=D\circ\prj_1+
\Phi_1\circ\prj_2$ and $q_\Psi=D\circ\prj_1+\Psi_1\circ\prj_2$, is
commutative. Now $\chi^{-1}(E_0\oplus C) =E_0\oplus C$ and the
conclusion follows from the definition of going down.
\qed

We observe also that since $q^{-1}(\bar C)=E_0\oplus C\subset 
E_0\oplus F_1$ and $q|_{0\oplus F_1}=\Phi_1$, $\Phi_1^{-1}(\bar C)=C$.
In other words:
\begin{prop}\label{left_inverse}
  Whenever going down is applicable to $\Phi_\bullet:F_\bullet
  \rightarrow E_\bullet$ then $\Phi_\bullet^!$ is a left inverse to
  $(\Phi_\bullet)_!$, that is
  \[
    \Phi_\bullet^!(\Phi_\bullet)_!(C)\ =\ C
  \]
  for any $F_0$-cone $C\subset F_1$.
\qed
\end{prop}
Note that $\Phi_\bullet^!$ is generally not right-inverse to
$(\Phi_\bullet)_!$. For example consider $\Phi_\bullet=(\Id,\iota):
(F_0\to F_1)\to (F_0\to F_1\oplus N)$ for any linear space $N$ over
$X$, $\iota: F_1\to F_1\oplus N$ the inclusion of the first factor and
$F_0$ acting trivially on $N$. Then for an $F_0$-cone of the form
$C\oplus N$ it holds $(\Phi_\bullet)_!\Phi_\bullet^!(C\oplus N) =
C\oplus 0$. Compare however Proposition~\ref{1-1-corr}.

Going down is functorial:
\begin{prop}
  Let $\Psi_\bullet:G_\bullet\rightarrow F_\bullet$, $\Phi_\bullet:
  F_\bullet\rightarrow E_\bullet$ be commutative squares of linear
  spaces to which going down is applicable, and let $C\subset G_1$ be a
  $G_0$-cone. Then
  \[
    (\Phi_\bullet\circ\Psi_\bullet)_!(C)
	\ =\ (\Phi_\bullet)_! (\Psi_\bullet)_!(C)\, .
  \]
\end{prop}
\pf
Consider the following diagram of linear spaces and cones:
\[
\begin{array}{ccccccccc}
&&0&&0\\[8pt]
&&\diagl{}&&\diagl{}\\[14pt]
&&G_0&\stackrel{\Id_{G_0}}{\llongrightarrow}&G_0\\[5pt]
&&\diagl{
\bigg(\!\!\!\begin{array}{c}
\scriptstyle\Psi_0\\ \scriptstyle\Id_{G_0}
\end{array}\!\!\!\bigg)}
&&\diagr{
\left(\!\!\!\begin{array}{c}
\scriptstyle\Psi_0\\ \scriptstyle 0\\ \scriptstyle -D_G
\end{array}\!\!\!\right)}\\[-30pt]
0&\llongrightarrow&F_0\oplus G_0
  &\stackrel{\left(\!\!\!\begin{array}{cc}
  \scriptstyle\Id_{F_0}& \scriptstyle 0\\
  \scriptstyle -\Phi_0& \scriptstyle \Phi_0\circ\Psi_0\\
  \scriptstyle 0& \scriptstyle -D_G
  \end{array}\!\!\!\right)}{\llongrightarrow}&F_0\oplus E_0\oplus C
  &\stackrel{(\Phi_1\circ D_F, D_E, \Phi_1\circ\Psi_1)}{\llongrightarrow}
  &(\Phi_\bullet\circ \Psi_\bullet)_! C&\llongrightarrow &0\\[8pt]
  &&\diagl{(-\Id_{F_0},\Psi_0)}&&\diagl{ \left(\!\!\!\begin{array}{ccc}
  \scriptstyle 0&\!\!\scriptstyle \Id_{E_0}&\!\! \scriptstyle 0\\
  \scriptstyle D_F&\!\!\scriptstyle 0&\!\! \scriptstyle \Psi_1
  \end{array}\!\!\!\right)}&&\diagl{\Id_{E_1}}\\[12pt]
0&\llongrightarrow&F_0&\llongrightarrow&E_0\oplus (\Psi_\bullet)_! C
  &\stackrel{(D_E,\Phi_1)}{\llongrightarrow}
  &(\Phi_\bullet)_!(\Psi_\bullet)_! C
  &\llongrightarrow &0\\[-3pt]
&&&\left(\!\!\!\begin{array}{c}
\scriptstyle\Phi_0\\ \scriptstyle -D_F
\end{array}\!\!\!\right)\\[-15pt]
&&\diagl{}&&\diagl{}\\[12pt]
&&0&&0
\end{array}
\]
We claim that rows and columns in this diagram are exact. The left
vertical sequence is trivially exact. The upper howizontal sequence is
exact by Proposition~\ref{ex_seq_cones} applied to $(\Phi_\bullet)_!$
and the cone $(\Psi_\bullet)_!(C)\subset F_1$. Exactness of the middle
vertical sequence follows by adding a trivial $E_0$-term to the
analogous sequence for $(\Psi_\bullet)_!$ and $C\subset G_1$.
For the remaining middle horizontal sequence
exactness of the enveloping sequence
\[
0\longrightarrow F_0\oplus G_0\longrightarrow
F_0\oplus E_0\oplus G_1 \longrightarrow E_1\longrightarrow 0 
\]
of linear spaces is easy to verify. Again by
Proposition~\ref{ex_seq_cones} the preimage $\tilde C\subset F_0\oplus
E_0\oplus G_1$ of $(\Phi_\bullet\circ\Psi_\bullet)_!(C)\subset E_1$
intersects $0\oplus E_0\oplus G_1$ in $0\oplus E_0\oplus C$, and it is
invariant under the action of $F_0$ on the first factor. Hence $\tilde
C= F_0\oplus E_0\oplus C$, proving exactness of the middle horizontal
sequence.

Now exactness of the lower horizontal sequence and of the middle
vertical sequence show that the preimage of
$(\Phi_\bullet)_!(\Psi_\bullet)_! C$ under the composition of
epimorphisms $F_0\oplus E_0\oplus G_1\to E_0\oplus F_1\to E_1$ from
the lower right square equals $F_0\oplus E_0\oplus C$. This is the
same as the preimage of $(\Phi_\bullet\circ \Psi_\bullet)_!C$.
Therefore $(\Phi_\bullet)_!(\Psi_\bullet)_! C$ and
$(\Phi_\bullet\circ\Psi_\bullet)_! C$ are the same cones in $E_1$.
\qed


\subsection{The case of quasi-isomorphisms}
By definition a morphism $\Phi_\bullet$ of two-term complexes is a
quasi-isomorphism if $H^i(\Phi_\bullet)$ is an isomorphism for
$i=0,1$. This is equivalent to requiring that $\Phi_\bullet$ viewed
as a commutative square is cartesian and cocartesian, see
Lemma~\ref{ass_exact_sequence}. Going up and down behaves well with
respect to quasi-isomorphisms:
\begin{prop}\label{1-1-corr}
  Let $\Phi_\bullet:F_\bullet\rightarrow E_\bullet$ be a 
  quasi-isomorphism of two-term complexes of linear spaces with
  $F_0$ locally free. Then going up and down induces a functorial
  one-to-one correspondence between $F_0$-cones $C\subset F_1$
  and $E_0$-cones $\bar C\subset E_1$.
\end{prop}
\pf
In view of Proposition~\ref{left_inverse} it remains to show that if 
$\bar C\subset E_1$ is an $E_0$-cone then $\bar C=(\Phi_\bullet)_!
\Phi_\bullet^!(\bar C)$. This is a local problem. We may thus assume
that there exists a local splitting $\sigma: E_0\oplus F_1\to F_0$ of
the exact sequence
\[
  0\longrightarrow F_0\longrightarrow E_0\oplus F_1
  \stackrel{q}{\longrightarrow}E_1\longrightarrow0,\quad
  q=D\circ \prj_1+\Phi_1\circ\prj_2
\]
from Lemma~\ref{ass_exact_sequence}.  Then $\chi=(\sigma,q):
E_0\oplus F_1\to F_0\oplus E_1$ is an isomorphism mapping the
diagonal $F_0$-action on $E_0\oplus F_1$ to the action on the first
factor of $F_0\oplus E_1$. Since $\sigma$ is a splitting, $\chi$
induces an isomorphism $\ker(\sigma)\to E_1$. Therefore
\[
\chi\big(\ker(\sigma)\cap q^{-1}(\bar C)\big)
=\chi(q^{-1}(\bar C))\cap (0\oplus E_1)=0\oplus \bar C.
\]
But $\chi(q^{-1}(\bar C))$ is an $F_0$-cone, and hence
Proposition~\ref{quot_cone} implies $\chi(q^{-1}(\bar C)) = F_0\oplus
\bar C$. By definition this says $\bar C=(\Phi_\bullet)_!
\Phi_\bullet^!(\bar C)$.
\qed

Using the nice behavior under quasi-isomorphisms we may now 
define {\em going down for morphisms in the derived category}
$D(\Coh(X))$ of the category of coherent sheaves. In the language of
linear spaces a morphism of two-term complexes in the derived
category $\Phi_\bullet: F_\bullet\rightarrow E_\bullet$ consists of
\begin{enumerate}
\item
  another two-term complex $G_\bullet$
\item
  a quasi-isomorphism $\Theta_\bullet:E_\bullet\rightarrow G_\bullet$
\item
  and a morphism $\Psi_\bullet:F_\bullet\rightarrow G_\bullet$.
\end{enumerate}
Two morphisms defined by tuples $(G_\bullet,\Theta_\bullet,
\Psi_\bullet)$ and $(G'_\bullet,\Theta'_\bullet,\Psi'_\bullet)$ are
considered equivalent if there exists a two-term complex ${\tilde
G}_\bullet$ and quasi-isomorphisms $\Lambda_\bullet:
G_\bullet\rightarrow{\tilde G}_\bullet$, $\Lambda'_\bullet:
G'_\bullet\rightarrow{\tilde G}_\bullet$ making the following diagram
commutative up to homotopy:
\[
\begin{array}{ccccc}
  &&\ G_\bullet\ &&\\[7pt]
  &\lbrtdiagl{\Psi_\bullet}&\diagr{\raisebox{-5pt}{$\scriptstyle
  \Lambda_\bullet$}}&\rbltdiagr{\Theta_\bullet}&\\[15pt]
  \hspace{5cm}F_\bullet\ &&{\tilde G}_\bullet&&\ E_\bullet
  \hspace{5cm}(*)\\[7pt]
  &\ltrbdiagl{\Psi'_\bullet}&\updiagr{\raisebox{2pt}{$\scriptstyle
  \Lambda'_\bullet$}}&\rtlbdiagr{\Theta'_\bullet}&\\[15pt]
  &&G'_\bullet&&
\end{array}
\]
\begin{defi}\label{go_down_deriv}\rm
  (going down in the derived category) Let $\Phi_\bullet:F_\bullet
  \rightarrow E_\bullet$ be a morphism of two-term complexes of
  linear spaces in the derived category, inducing an isomorphism on
  $H^0$ and a closed embedding on $H^1$. Moreover, we require
  $F_0$ to be locally free. When these assumptions are 
  satisfied we say that {\em going down is applicable to $\Phi_\bullet$}.
  In this case the {\em push-forward} of an $F_0$-cone $C\subset F_1$
  is defined to be the $E_0$-cone
  \[
    (\Phi_\bullet)_!(C):=(\Theta_\bullet)^!(\Psi_\bullet)_!(C)\subset 
    E_1\, ,
  \]
  whenever $(G_\bullet,\Psi_\bullet,\Theta_\bullet)$ is a 
  representative of $\Phi_\bullet$.
\qed
\end{defi}
Using the previous results it is easy to check that this is well-defined.
\begin{rem}\rm\label{coh_cones}
  One might wonder if there exists  ``going down'' when being only
  given maps on the level of cohomology. There are three remarks I
  want to make on this.
  \begin{enumerate}
  \item
    A map in cohomology is considerably weaker than a map of 
    complexes, even for two-term complexes. For instance, let
    $E_\bullet=[E_0\stackrel{D}{\rightarrow}E_1]$ be a non-split
    epimorphism of linear spaces and $K=\ker D$. Then
    $H^\bullet(E_\bullet) =H^\bullet([K_\bullet\rightarrow0])$, but the
    identity map in cohomology is not induced by a morphism of
    complexes.
  \item
    There is going down for ``cones coming from cohomology'': 
    By such cones we mean cones of the form $C=p^{-1}(\bar C)$ for
    some $\bar C\subset H^1(E_\bullet)$, $p:E_1 \rightarrow H^1
    (E_\bullet)$ the cokernel of $E_\bullet$. Namely, if $\rho: H^1(
    E_\bullet)\rightarrow H^1(F_\bullet)$ is a closed embedding and
    $p': F_1 \rightarrow H^1(F_\bullet)$ is the cokernel of $F_\bullet$,
    one may set $\rho_!(C) :={p'}^{-1}(\rho(\bar C))$. In case
    $\rho=H^1(\Phi_\bullet)$ with $\Phi_\bullet:E_\bullet\rightarrow
    F_\bullet$ a morphism to which going down is applicable, then
    $\rho_!(C)$ obviously coincides with $(\Phi_\bullet)_!(C)$.
  \item
	Not every $E_0$-cone in $E_1$ comes from cohomology. As a simple
	example with $C\neq E_1$ but $H^1(E_\bullet)=0$ take
	$X=\Aff^1_k=\Spec k[T]$, $E_0=E_1=L(\calo_X) =\Spec k[T,U]$,
	$D:E_0\rightarrow E_1$ corresponding to the homomorphism of
	$k[T]$-algebras sending $U$ to $TU$, and $C=V(TU)\subset E_1$, the
	linear space corresponding to the structure sheaf of the origin.

    See also Remark~\ref{C_X|M_descend?} for another, less artificial 
    example.
  \qed
  \end{enumerate}
\end{rem}


\section{Global normal cones}
If $\iota:X\hookrightarrow M$ is a closed embedding of algebraic 
$k$-schemes the normal cone $C_{X|M}\subset N_{X|M}$ is a
$TM|_X$-cone (Example~\ref{C_Z|X}). With nonsingular $M$ theses
normal cones are essentially unique, namely up to vector bundle
factors. In fact, if $\iota':X\hookrightarrow M'$ is another such
embedding we may consider the diagonal $(\iota,\iota'):X
\hookrightarrow M\times M'$ to reduce to the case where $\iota=
\pi\circ\iota'$, $\pi:M'\rightarrow M$ a smooth morphism. But
then there is an exact sequence of cones
\begin{eqnarray}
  0\longrightarrow(\iota')^*T_{M'|M}\longrightarrow C_{X|M'}
  \longrightarrow C_{X|M}\longrightarrow0\, .\label{eqn1}
\end{eqnarray}
Based on this observation Behrend and Fantechi show that to any 
$X$ there is associated a cone stack (a certain Artin stack) over $X$ of
pure relative dimension zero, the {\em intrinsic normal cone}
$\calc_X$. Locally, $\calc_X$ is nothing but the stack-theoretic
quotient $C_{X|M}/T_M|_X$, and the above exact sequence of cones is
responsible for the fact that these quotients glue.

One essential insight of Behrend and Fantechi is that one can retrieve 
an actual cone over $X$ by giving a morphism $\ph^\bullet:
\calf^\bullet\rightarrow\call_X^\bullet$ in the derived category
inducing an isomorphism in $H^0$ and an epimorphism in $H^{-1}$, and
such that $\calf^\bullet=[\calf^{-1}\rightarrow\calf^0]$ is a two-term
complex of locally free sheaves. Here $\call_X^\bullet$ is the
cotangent complex of $X$. (In the language of
\cite{behrendfantechi}, $\ph^\bullet$ is a ``global resolution'' of a
``perfect obstruction theory'' for $X$.)

The cotangent complex is a complicated and largely mysterious object
canonically associated to any scheme, or even ringed topos
\cite{illusie}. However, here we will work exclusively with the
truncated complex  $\tau_{\ge -1}\call_X^\bullet$. This is simply an
object of the derived category that has the following explicit local
description: If $U\subset X$ is an open subscheme and
$U\hookrightarrow M$ is a closed embedding into a smooth scheme $M$
then
\[
\tau_{\ge -1} \call_X^\bullet= [\cali/\cali^2\to \Omega_M|_U],
\]
where the complex on the right hand side has entries at $-1$ and $0$
(this follows from the exact triangle for the cotangent complex, see
below). In particular, if $X$ is globally embedded into a smooth
scheme we can avoid the cotangent complex at all.

Using our study of going up and down for $E$-cones we will see that 
the object needed is the following.
\begin{defi}\label{global_normal_space}\rm
  A {\em global normal space} for $X$ is a morphism $\ph^\bullet:
  \calf^\bullet = [\calf^{-1}\rightarrow\calf^0] \rightarrow
  \tau_{\ge -1}\call_X^\bullet$ in the derived category with
  $\calf^0$ locally free and inducing an isomorphism in $H^0$ and an
  epimorphism in $H^{-1}$.
\qed
\end{defi}

Given a global normal space $\Phi_\bullet:\tau_{\le1}
(L_X)_\bullet\rightarrow  F_\bullet$, now written in terms of linear
spaces $\tau_{\le 1}(L_X)_\bullet=L(\tau_{\ge -1}\call_X^\bullet)$
etc., we may construct a cone $C=C(\Phi_\bullet)\subset F_0$ as
follows: Let $U\subset X$ be an open set embedded into a nonsingular
$M$,  $\iota: U\hookrightarrow M$. The exact triangle of relative
cotangent complexes associated to $U\rightarrow M\rightarrow\Spec k$
yields a morphism in the derived category
\[
  \Lambda_\bullet:[T_M|_U\rightarrow N_{U|M}]\longrightarrow
  \tau_{\le 1}(L_U)_\bullet
\]
that induces isomorphisms in $H^i$, $i=0,1$. As $T_M|_U$
is a vector bundle the composition  $\Phi_\bullet|_U\circ
\Lambda_\bullet$ fulfills the assumption of
Definition~\ref{go_down_deriv}. We may thus define
\begin{eqnarray}
  C|_U\ :=\ \big(\Phi_\bullet|_U\circ\Lambda_\bullet\big)_!
  (C_{U|M})\, .\label{eqn2}
\end{eqnarray}
It remains to show
\begin{lemma}\label{independence}
  $(\Phi_\bullet|_U\circ\Lambda_\bullet)_!(C_{U|M})\subset F_1|_U$ is 
  independent of choices.
\end{lemma}
\pf
It suffices to treat the case of another embedding $\iota':U
\rightarrow M'$ s.th.\ $\iota=\pi\circ\iota'$ for some smooth 
morphism $\pi:M'\rightarrow M$, see above. We have a commutative
diagram with exact rows and columns
\[
\begin{array}{cccccccccccl}
  &&&&0&&0\\[8pt]
  &&&&\diagl{}&&\diagl{}\\[12pt]
  &&&&\makebox[0pt]{${\iota'}^*T_{M'|M}$}&=\!=\!=&
    \makebox[0pt]{${\iota'}^*T_{M'|M}$}\\[8pt]
  &&&&\diagl{}&&\diagl{}\\[12pt]
    0&\llongrightarrow&T_U&\llongrightarrow&{\iota'}^*T_{M'}&
    \llongrightarrow&N_{U|M'}
    &\llongrightarrow&T_1(U)&\llongrightarrow&0\\[8pt]
  &&\diagl{}&&\diagl{}&&\diagr{}&&\diagr{}\\[12pt]
  0&\llongrightarrow&T_U&\llongrightarrow&\iota^*T_M&
    \llongrightarrow&N_{U|M}
    &\llongrightarrow&T_1(U)&\llongrightarrow&0&\!\!\!.\\[8pt]
  &&&&\diagl{}&&\diagl{}\\[12pt]
  &&&&0&&0
\end{array}
\]
Here $T_1(U)$ is the linear space associated to the first higher cotangent 
sheaf of $U$. This shows that $D \pi$ induces a quasi-isomorphism
$\Theta_\bullet: [T_{M'}|_U \rightarrow N_{U|M'}] \rightarrow
[T_M|_U\rightarrow N_{U|M}]$ with $\Lambda'_\bullet=
\Lambda_\bullet \circ\Theta_\bullet$. Moreover, from the exact
sequence of cones (\ref{eqn1})
\[
  C_{U|M'}\ =\ \Theta_1^{-1}(C_{U|M})\ =\ \Theta_\bullet^!(C_{U|M})\, .
\]
Thus by Proposition~\ref{1-1-corr} we conclude
\[
  (\Lambda'_\bullet)_!(C_{U|M'})\ =\ 
  (\Lambda_\bullet)_!(\Theta_\bullet)_!\Theta_\bullet^!(C_{U|M})
  \ =\ (\Lambda_\bullet)_!(C_{U|M})\, .
\]
\qed

Here is the first of the two main results of this paper.
\begin{theorem}\label{global_normal_cone}
  Let $X$ be an algebraic $k$-scheme. To any global normal space
  $\Phi_\bullet: \tau_{\le 1}(L_X)_\bullet\rightarrow F_\bullet$ for $X$ is 
  associated an $F_0$-cone $C(\Phi_\bullet)\subset F_1$, locally of the
  form (\ref{eqn2}), and of pure dimension equal to $\rk F_0$.
\end{theorem}
\pf
It remains to check the statement on the dimension. Locally, we 
may choose a representation of $\Phi_\bullet\circ \Lambda_\bullet:
[T_M|_U\rightarrow N_{U|M}] \rightarrow F_\bullet$, where
$\iota:U\hookrightarrow M$ is a closed embedding of an open
$U\subset X$, by $(G_\bullet,\Theta_\bullet,\Psi_\bullet)$ with
\begin{itemize}
\item
  $\Psi_\bullet:[T_M|_U\rightarrow N_{U|M}]\rightarrow G_\bullet$
\item
  $\Theta_\bullet:F_\bullet\rightarrow G_\bullet$ a quasi-isomorphism
\item
  $G_\bullet=[G_0\rightarrow G_1]$ with $G_0$ a vector bundle (!).
\end{itemize}
We get two exact sequences of cones
(Proposition~\ref{ex_seq_cones})
\[
\begin{array}{rcccccccl}
  0&\longrightarrow&T_M|_U&\longrightarrow&G_0\oplus C_{U|M}
  &\longrightarrow&(\Psi_\bullet)_!C_{U|M}&\longrightarrow&0\\[1ex]
  0&\longrightarrow&F_0&\longrightarrow&G_0\oplus 
  C(\Phi_\bullet)|_U
  &\longrightarrow&(\Psi_\bullet)_!C_{U|M}&\longrightarrow&0\, .
\end{array}
\]
The first one shows that $(\Psi_\bullet)_!C_{U|M}$ is pure dimensional 
of dimension equal to $\rk G_0$, and then by the second one
$C(\Phi_\bullet)$ is pure $(\rk F_0)$-dimensional.
\qed

\begin{defi}\rm
  $C(\Phi_\bullet)$ is called the {\em global normal cone} associated 
  to the global normal space $\Phi_\bullet$.
\qed
\end{defi}
\begin{rem}\rm\label{C_X|M_descend?}
1)\ \  The picture would be especially simple if for any closed
embedding $\iota:X \hookrightarrow M$ into a nonsingular $M$,
$C_{X|M}$ came from a cone in the intrinsically defined linear space
$T_1(X)$. This is, however, generally wrong:

Consider the fat point $X=\Spec R$, $R=k[X,Y]/(X^2,XY,Y^2)$, with its 
embedding into $M=\Aff_k^2=\Spec k[X,Y]$. Letting $A$, $B$, $C$
correspond to the generators $X^2$, $XY$, $Y^2$ of the ideal,
$C_{X|M}=\Spec R[A,B,C]/ (B^2-AC,XC-YB,XB-YA)$. $T_1(X)$ is the linear
space corresponding to the kernel $I_T$ of
\begin{eqnarray*}
  &R[A,B,C]/(XC-YB,XB-YA)\ \longrightarrow\ R[dX,dY],\\
  &A\mapsto 2XdX,\ B\longmapsto YdX+XdY,\ C\mapsto 2YdY\, ,
\end{eqnarray*}
which is $(XA,YA,XB,YB,XC,YC)=(X,Y)\cdot(A,B,C)$. A cone in
$N_{X|M}=\Spec R[A,B,C]/(XC-YB,XB-YA)$ comes from  $T_1(X)$ iff its
ideal is generated by polynomials in $XA$, $YA$, $XB$, $YB$, $XC$,
$YC$. This is not the case for $B^2-AC$, and so $C_{X|M}$ does not
come from a cone in $T_1(X)$.\\[1ex]
2)\ \ By uniqueness of minimal free resolutions of modules
over a local ring (see e.g.\ \cite[Thm.20.2]{eisenbud}) it is not hard to
show that for any, not necessarily closed point $x\in X$ there is a
{\em minimal} germ of global normal spaces at $x$. This is
constructed by embedding an \'etale neighborhood of $x$ in $X$ into a
smooth $k$-scheme $M$ of dimension ${\rm embdim}_x X= \dim_{k(x)}
\Omega_x\otimes k(x)$. We assume $k$ perfect here to assure that
regular $k$-schemes are smooth over $k$. A germ of global
normal space at $x$ can then be defined by selecting a minimal set of
generators for the ideal defining $X\hookrightarrow M$. This
germ of global normal space is minimal in the sense that any
other germ of global normal space at $x$ can be obtained by adding
trivial factors.

As a consequence, the germ at $x$ of any global normal cone is
isomorphic to $C_{X|M}$ plus a vector bundle factor. Morally speaking,
the ``nonlinear parts'' of global normal cones are locally unique.
\qed
\end{rem}


\section{Virtual fundamental class and Fulton's canonical class}
\subsection{Virtual fundamental classes}
If $X$ is an algebraic $k$-scheme and $\Phi_\bullet:\tau_{\le 1}
(L_X)_\bullet \rightarrow F_\bullet$ is a global normal space for $X$
with also {\em $F_1$ a vector bundle} we speak of a {\em free} global
normal space of {\em rank} $\rk(\Phi_\bullet)=\rk F_0-\rk F_1$. We
may then intersect the zero section $s:X\rightarrow F_1$ of $F_1$
with the global normal cone $C(\Phi_\bullet)\subset F_1$ to produce a
class on $X$.
\begin{defi}\rm
  Let $\Phi_\bullet: \tau_{\le 1} (L_X)_\bullet\rightarrow F_\bullet$
  be a free  global normal space. The Chow class
  \[
    [X,\Phi_\bullet]\ :=\ s^![C(\Phi_\bullet)]\in A_{\rk(\Phi_\bullet)}(X)
  \]
  is called {\em virtual fundamental class} of $X$ with respect to
  $\Phi_\bullet$.
\qed
\end{defi}
Note that $[X,\Phi_\bullet]$ contains as much information as 
$[C(\Phi_\bullet)]\in A_{\rk(F_0)}(F_1)$, for
\[
  [C(\Phi_\bullet)]\ =\ p^!s^![C(\Phi_\bullet)]\ =\ p^![X,\Phi_\bullet]\, ,
\]
where $p:F_1\rightarrow X$ is the projection.

One of the most important property of such classes is their 
compatibility with specializations. In the application to the
construction of invariants from moduli spaces associated to a
projective manifold V, say (as in Gromov-Witten or Donaldson-theory),
this property implies invariance under smooth deformations of $V$.
There are two versions of the specialization theorem, one involving
global normal spaces of the total space of a family, and the other
working with {\em relative} global normal spaces (that is, a morphism
$\Phi_\bullet: \tau_{\le 1}(L_{\calx|S})_\bullet \rightarrow
F_\bullet$, where $\calx\rightarrow S$ is the family of algebraic
$k$-schemes under consideration). We do not have anything to add to
the presentation in \cite[Proposition~5.10 and
Proposition~7.2]{behrendfantechi}, which translates  almost literally
into our language.

Before turning to an explicit formula for the computation of 
$[X,\Phi_\bullet]$ in hopefully more accessible terms, we want to add
the following point of view: For any pure-dimensional cone $C$ in a
vector bundle $F$ there is a formula for the intersection with the zero
locus in terms of the Segre class of $C$ and the total Chern class of
$F$ \cite[Expl.\ 4.1.8]{fulton}. Applied to $C(\Phi_\bullet)\subset
F_1$ it says
\[
  [X,\Phi_\bullet]\ =\ \{c(F_1)\cap s(C(\Phi_\bullet))
  \}_{\rk(\Phi_\bullet)}\, ,
\]
where $\{\,.\,\}_d: A_*(X)\rightarrow A_d(X)$ denotes the projection
to the $d$-dimensional part. Now for any $r>0$, the image of
$[C(\Phi_\bullet)]$ under the  monomorphism $\iota_r: F_1
\hookrightarrow F_1\oplus\Aff_X^r$ becomes rationally trivial, while
$c(F_1\oplus \Aff^r_X) = c(F_1)$. Thus letting ${\tilde
s}^r:X\hookrightarrow F_1\oplus\Aff_X^r$ be the zero section we see
\[
  \{c(F_1)\cap s(C(\Phi_\bullet))\}_{\rk(\Phi_\bullet)-r}\ =\ 
  ({\tilde s}^r)^!(\iota_r)_*[C(\Phi_\bullet)]\ =\ 0\, .
\]
This teaches us two things: First, if $F_1$ splits off a trivial factor
$F_1={\bar F}_1\oplus\Aff^1_X$ with $\img\Phi_1 \subset{\bar F}_1$
then $[X,\Phi_\bullet]=0$. So the result is trivial if $F_1$ is not chosen
small enough. And second, if $[X,\Phi_\bullet]\neq 0$ then
$\rk(\Phi_\bullet)$ is the smallest number $d$ such that
\[
  \{c(F_1)\cap s(C(\Phi_\bullet))\}_d\ \neq\ 0\, .
\]
\begin{expl}\rm
  Let $X$ be smooth of dimension $n$. Then the cotangent complex of 
  $X$ is exact at $\call_X^{-1}$. So the natural morphism
  $\Phi_\bullet: \tau_{\le 1} (L_X)_\bullet\rightarrow
  [T_X\rightarrow O]$, $O=X$ the trivial linear space over $X$, is an
  isomorphism in $H^0$ and $H^1$. Then $C(\Phi_\bullet)=X=O$ and
  $c(F_1)\cap s(C(\Phi_\bullet)) =s(C(\Phi_\bullet))=[X]$ has
  vanishing components in dimensions  smaller than $n$.
\qed
\end{expl}

\subsection{Fulton's canonical class}
If an algebraic $k$-scheme $X$ is globally embeddable into a 
smooth $k$-scheme $M$ (e.g.\ $X$ quasi-projective) then
\[
  c_F(X)\ :=\ c(T_M|_X)\cap s(C_{X|M})\, \in A_*(X)
\]
is a Chow-class on $X$ that is independent of the choice of 
embedding \cite[Expl.\ 4.2.6]{fulton}.
\begin{defi}\label{fultons_class}\rm
  $c_F(X)$ is called {\em Fulton's canonical class}.
\qed
\end{defi}
Note that if $X$ is smooth one may choose $X=M$ and so 
$c_F(X)=c(T_X)\cap[X]$. For comparison of $c_F(X)$ with Mather's
and MacPherson's Chern classes see \cite{aluffi}.

Given a (not necessarily free) global normal space $\Phi_\bullet:
\tau_{\le 1} (L_X)_\bullet\rightarrow  F_\bullet$ for $X$, $c_F(X)$
can also be expressed as follows:
\begin{prop}\label{Fultons_class}
  Let $\Phi_\bullet:\tau_{\le 1} (L_X)_\bullet\rightarrow F_\bullet$
  be a global  normal space for a quasi-projective $X$. Then
  \[
    c_F(X)\ =\ c(F_0)\cap s(C(\Phi_\bullet))\, .
  \]
\end{prop}
\pf
By quasi-projectivity there exists a {\rm global} closed embedding
$\iota:X\hookrightarrow M$ of $X$ into a smooth $M$. This yields a
globally defined morphism in the derived category $\Lambda_\bullet:
[T_M|_X\rightarrow N_{X|M}]\rightarrow \tau_{\le 1} (L_X)_\bullet$.
Also by quasi-projectivity any sheaf is the quotient of a locally
free sheaf. Hence there is a global representative $(G_\bullet,
\Theta_\bullet, \Psi_\bullet)$ of $\Phi_\bullet\circ \Lambda_\bullet$
in the construction of Theorem~\ref{global_normal_cone}, that is,
with $G_0$ a vector bundle, $\Theta_\bullet:F_\bullet \rightarrow
G_\bullet$ a quasi-isomorphism, $\Psi_\bullet: [T_M|_X\rightarrow
N_{X|M}]\rightarrow G_\bullet$. We get two exact sequences of cones
with vector bundle kernels (see  Proposition~\ref{ex_seq_cones})
\[
\begin{array}{rcccccccl}
  0&\longrightarrow&T_M|_X&\longrightarrow&
    G_0\oplus C_{X|M}&\longrightarrow&
    (\Psi_\bullet)_!(C_{X|M})&\longrightarrow&0\\[2ex]
  0&\longrightarrow&F_0&\longrightarrow&G_0\oplus 
  C(\Phi_\bullet)&\longrightarrow&(\Psi_\bullet)_!(C_{X|M})
  &\longrightarrow&0\, ,
\end{array}
\]
which by the multiplicativity of Segre classes in exact 
sequences of cones with vector bundle kernels imply
\[
  c(T_M|_X)\cap s(C_{X|M})\ =\ c(G_0)\cap s\Big((\Psi_\bullet)_!
  (C_{X|M})\Big)\ =\ c(F_0)\cap s(C(\Phi_\bullet))\, .
\]
\vspace{-6ex}

\qed
\begin{rem}\rm
  If $X$ is any algebraic $k$-scheme with global normal spaces one 
  could take the right-hand side of the formula in the proposition as
  definition for a generalization of Fulton's canonical class on
  projective schemes. However, I was not able to prove independence
  of this class from the choice of $\Phi_\bullet$. And in case $X$ is
  not quasi-projective but embeddable into a smooth scheme, in the
  construction of Theorem~\ref{global_normal_cone} we might not
  be able to choose $G_0$ locally free. Then the coincidence of this
  class with $c_F(X)$ is not clear either. The problem is that on one
  hand the globally defined complex linking two global normal spaces
  $\Phi_\bullet: \tau_{\le 1}(L_X)_\bullet\rightarrow F_\bullet$,
  $\Phi'_\bullet: \tau_{\le 1} (L_X)_\bullet \rightarrow F'_\bullet$
  is the cotangent complex, which need not be globally representable
  by a complex $L_\bullet$ with $L_0$ a vector bundle, while on the
  other hand Segre classes do not behave well in exact sequences
  unless the kernels are vector bundles.
\qed
\end{rem}

We are now ready to deduce the announced formula for the 
virtual fundamental class.
\begin{theorem}\label{closed_formula}
  Let $X$ be a projective $k$-scheme and $\Phi_\bullet: \tau_{\le 1}
  (L_X)_\bullet\rightarrow F_\bullet$ a free global normal space for
  $X$ of constant rank $d$. Then
  \[
    [X,\Phi_\bullet]\ =\ \Big\{c(\ind F_\bullet)^{-1}\cap c_F(X)
    \Big\}_d\, ,
  \]
  where  $\ind F_\bullet$ is the virtual bundle $F_0-F_1\in K^0(X)$.
\end{theorem}
\pf
As remarked at the end of the last subsection the virtual 
fundamental class can be computed by the formula
\[
  s^![C(\Phi_\bullet)]\ =\ \Big\{c(F_1)\cap 
  s(C(\Phi_\bullet))\Big\}_d\, .
\]
Now just insert $c(F_0)^{-1}\cup c(F_0)$ and use 
Proposition~\ref{Fultons_class}
\qed
\begin{rem}\rm
1)\ \ 
  This formula enlightens the dependence of virtual fundamental
  classes  on the choice of global normal spaces: Interestingly,
  $[X,\Phi_\bullet]$ depends only on the index bundle of $F_\bullet$
  rather than on any of the finer data used to construct
  $C(\Phi_\bullet)$. But note also that for another choice
  $\Phi'_\bullet: \tau_{\le 1} (L_X)_\bullet\rightarrow F'_\bullet$
  of global normal space, $[X,\Phi'_\bullet]$ can not in general be
  computed from $[X,\Phi_\bullet]$ and $\ind F_\bullet$, $\ind
  F'_\bullet$ alone.\\[1ex]
2)\ \ 
  One can take this formula as {\em definition} of the virtual 
  fundamental class of $X$ without knowing anything about the more
  sophisticated theory of global normal cones in the non-projective
  case. This was the point of view of the author in summer 1995 in an
  attempt to define GW-invariants in algebraic geometry, when I
  observed it from formal considerations. Unfortunately, I was not
  aware of Vistoli's rational equivalence \cite{vistoli}, from which
  the crucial independence of the invariants under smooth
  deformations can be derived. I learned also that the same formula
  has independently been discovered by Brussee for complex spaces
  constructed as zero locus of holomorphic Fredholm sections of
  holomorphic Banach bundles over complex Banach manifolds, as
  occurring for example in Seiberg-Witten theory \cite{brussee} (the
  interpretation of Brussee's $c_*(X)$ as Fulton's canonical class is
  not quite clear, though).\\[1ex]
3)\ \ 
  At the beginning of Section~3 we mentioned Behrend and
  Fantechi's intrinsic normal cone $\calc_X$, which
  locally was the stack-theoretic quotient of $C_{X|M}$ by the action of
  $T_M|_X$ for some embedding $X\hookrightarrow M$ into a smooth
  space.
  
  Now if $X$ is globally embedded into a smooth scheme $M$, $\calc_X$ 
  is globally the quotient of $C_{X|M}$ by $T_M|_X$. Hence in view of 
  the multiplicative behavior of Segre classes in exact sequences of
  cones with vector bundle kernels, $c_F(X)$ could with some right
  considered as {\em Segre class of $\calc_X$}.

  Conversely, if there was a theory of Segre classes for cone stacks,
  the Segre class of $\calc_X$ would generalize Fulton's canonical class
  to arbitrary algebraic $k$-schemes.
\qed
\end{rem}

\end{document}